\documentclass[11pt,reqno]{amsart}
\usepackage{amsbsy,amsfonts,amsmath,amssymb,amscd,amsthm}
\usepackage{mathrsfs,float}
\usepackage[mathcal]{euscript}
\usepackage{enumitem,graphicx,epstopdf}
\usepackage[abs]{overpic}
\usepackage[a4paper,text={6.1in,8in},centering]{geometry}
\usepackage{pxfonts}
\usepackage[colorlinks=true,linkcolor=blue,citecolor=green]{hyperref}
\usepackage{srcltx}
\usepackage[margin=20pt,font=small,labelfont=bf,textfont=it]{caption}
\usepackage{graphicx}
\usepackage{subcaption}
\usepackage{pinlabel}
\usepackage[latin1]{inputenc}

\hypersetup{linktocpage}

\small\normalsize
\theoremstyle{plain}
\newtheorem{mainthm}{Theorem}
\newtheorem{maincor}{Corollary}
\newtheorem{thm}{Theorem}[section]
\newtheorem{lem}[thm]{Lemma}

\newtheorem{prop}[thm]{Proposition}
\newtheorem*{thm*}{Theorem}

\newtheorem{defi}[thm]{Definition}

\theoremstyle{definition}

\DeclareMathOperator{\uindex}{u-index}
\newcommand{\eqdef}{\stackrel{\scriptscriptstyle\rm def}{=}}

\setlist[enumerate,1]{label=(\arabic*)}
\setlist[enumerate,2]{label=(\alph*)}

\begin{document}

~\vspace{-1cm}
\title[]{Historic wandering domains
near cycles }

\dedicatory{Dedicated to Professor Lorenzo J.~D\'iaz on the
occasion of his $60$th birthday}

\author[Barrientos]{Pablo G.~Barrientos}
\address[Barrientos]{Instituto de Matem\'aticas e Estat\'istica, Universidade Federal Fluminense, Gragoata Campus, Rua Prof.\ Marcos Waldemar de Freitas Reis, S/n-Sao Domingos, Niteroi - RJ, 24210-201, Brazil}
\email{pgbarrientos@id.uff.br}

\begin{abstract}
We explain how to obtain non-trivial historic contractive
wandering domains for a dense set of diffeomorphisms in  certain
type of {$C^r$}-Newhouse domains of homoclinic tangencies in dimension
$d\geq 3$ {and $r\geq 1$. In particular, this gives for the first time a contribution to Takens' last problem in the $C^1$ topology and in dimension $d>2$.}
We show how these Newhouse domains can be obtained
arbitrarily close to diffeomorphisms exhibiting heterodimensional cycles {(in dimension $d=3$)}
or non-transverse equidimensional cycles {(in any dimension $d\ge 3$)} associated with periodic points with non-real
complex {leading} multipliers.
\end{abstract}

\maketitle \thispagestyle{empty}
\section{Introduction}

A \emph{non-trivial historic contractive wandering domain} for a
given map $f$ on a  $C^\infty$ Riemannian compact manifold $M$ is
a non-empty connected open set $D \subset M$ which satisfies the
following conditions:
\begin{itemize}
\item  $f^i(D)\cap f^j(D)=\emptyset$ for $i,j\geq 0$ with $i\not = j$,
\item  the union of the $\omega$-limit set for points  in $D$ is
not equal to a single orbit,
\item  the diameter of $f^i(D)$ converges to zero if $i\to \infty$,
\item the orbit of any point $x$ in $D$ has historic behavior,
i.e., the sequence of empirical measure
$\mu_n=\frac{1}{n}\sum_{i=0}^{n-1} \delta_{f^i(x)}$ does not
converge in the weak* topology.
%there exists a continuous map
%$\phi: M\to \mathbb{R}$ such that the sequence of Birkhoff average
%$\frac{1}{n} \sum_{i=0}^{n-1} \phi(f^i(x))$ does not converges.
\end{itemize}
Non-trivial contractive wandering domains were early observed by
Bohl and Denjoy (see \cite{Bohl,denjoy1932}) for $C^1$
diffeomorphisms on a circle. Following these results, similar
phenomena were observed for high dimensional examples {as well as one-dimensional maps in real analytic
category}~\cite{bonatti1994wandering, harrison1989denjoy, KNS17, knill1981c,  K10, Mc93, NS96, S85,L89b,L89a}. However, these domains are not
historic in the sense of the last condition above. The existence
of non-trivial historic contractive wandering domains
%in nonhyperbolic dynamics
were first studied by Colli and Vargas~\cite{CV01} for some
two-dimensional example which is made up of an affine thick
horseshoe with $C^2$-robust homoclinic tangencies. More recently
in~\cite{KS17} (see also \cite{berger2020emergence}) it was proved
that any two-dimensional diffeomorphism in any $C^2$-Newhouse
domain (open sets of $C^2$-diffeomorphisms with robust homoclinic
tangencies) is contained in the closure of diffeomorphisms having
non-trivial historic contractive wandering domains. In this paper,
we will explain how this result could be generalized to higher
dimensions {and in the $C^1$ topology} for certain class of Newhouse domains.
This provides {for the first time}  examples of smooth dynamical systems  {in dimension greater than two,
as well as in the $C^1$ topology, where it is not possible to get rid of historical
behavior by eliminating negligible sets of diffeomorphisms and of
initial conditions. This problem was raised by Ruelle in~\cite{Ruelle01} and similarly by Takens in~\cite{Takens08} being nowadays known  as \emph{Takens' last problem}.} Also answering~\cite[Question~2]{BKNRS20}, we give some conditions
ensuring that diffeomorphisms with certain types of
heterodimensional and equidimensional cycles have historic
contractive non-trivial wandering domains. %\vspace{-0.1cm}

%\section{Multidimensional historic non-trivial wandering domains}
%
%In this section we will explain the state of the art on Newhouse
%domains, Newhouse phenomenon and coexistence of infinitely many
%invariant circles. At the end, we explain how we can obtain more
%or less straight forward from the results in the literature
%multidimensional historic non-trivial wandering domains
%arbitrarily $C^r$-close to certain non-transverse equidimensional
%cycles and $C^1$-close to 3-dimensional heterodimensional cycles
%with complex eigenvalues.

\subsection{Wandering domains for Newhouse domains in higher dimensions}
Following~\cite{BD12}, we say that a $C^r$-open set $\mathcal{N}$
of diffeomorphisms is a \emph{$C^r$-Newhouse domain} if there
exists a dense set $\mathcal{D}$ in $\mathcal{N}$ such that every
$g\in \mathcal{D}$ has a homoclinic tangency associated with some
hyperbolic periodic saddle. Furthermore, if these homoclinic
tangencies satisfy a given property~$\mathcal{P}$, then we may
call it a $C^r$-Newhouse domain of homoclinic tangencies
satisfying~$\mathcal{P}$.  {See Definition~\ref{def-cycles} for a formal definition of homoclinic tangency (and heterodimensional cycle).}

The first example of a $C^r$-Newhouse
domain was obtained by Newhouse~\cite{New70} in any surface for
$r\geq 2$. Multidimensional
$C^r$-Newhouse domains for $r\geq 2$ was constructed by Palis and
Viana~\cite{PV94}, Romero~\cite{Ro95} and Godchenko, Shil'nikov
and Turaev in~\cite{GTS93} (see also~\cite{GST08}). Namely, from
these papers, %it follows that
$C^r$-Newhouse domains with $r\geq 2$
can be constructed in any manifold of dimension $d\geq 2$
arbitrarily $C^r$-close to any $C^r$-diffeomorphism having a
homoclinic tangency associated with a hyperbolic periodic point~$P$. This means that if $f$ is a $C^r$-diffeomorphism ($r\geq 2$)
with a homoclinic tangency associated with $P$, then $f\in
\overline{\mathcal{N}}$ where $\mathcal{N}$ is a $C^r$-Newhouse
domain of homoclinic tangencies. Moreover, these homoclinic tangencies are associated with periodic points
satisfying similar multiplier condition as $P$. The same result was
previously established by Newhouse in~\cite{New79} for surface
dynamics.

One of the first examples of $C^1$-Newhouse domains was
obtained in three-dimensional manifolds by Bonatti and D\'iaz
in~\cite{BD99} associated with homoclinic tangencies to periodic
points with complex eigenvalues and involving heterodimensional
cycles\footnote{Some comments on~\cite{BD99} are necessary. This
paper appeared before the stabilization theory of
heterodimensional cycles~\cite{BD08,BDK12}. The trick used by
Bonatti and D\'iaz was to consider a robust heterodimensional {cycle}
coming from~\cite{BD96} (where blenders were introduced) and
additionally they assumed that this cycle is homoclinically
related to another heterodimensional cycle with complex
eigenvalues. They showed that, by a $C^1$-perturbation, the
homoclinic classes of the involved periodic points are
$C^1$-robustly linked. This provides a $C^1$-open set where
densely there exist homoclinic tangencies associated with periodic
points with complex multipliers. That is, they construct a
$C^1$-Newhouse domain of homoclinic tangencies of periodic points
with complex multipliers. At that time, they did not know if this
open set corresponds with an open set of robust tangencies (i.e.,
where the homoclinic tangencies associated with a non-trivial
hyperbolic set persists under perturbations). But now,
from~\cite{BD12} we can conclude that, indeed, $C^1$-robust
tangencies associated with a blender-horseshoe appear in this
open set.}. Later, Asaoka in~\cite{Ao08} provides $C^1$-Newhouse
domains in any manifold of dimension $d\geq 3$. As Asaoka himself
mentioned in~\cite{Ao09}, his example is essentially the same that
Simon previously provided in~\cite{Si72}. Both examples are based
on normally hyperbolic non-trivial attractors. More recently,
again Bonatti and D\'iaz in~\cite{BD12} have constructed similar
examples of $C^1$-Newhouse domains but now associated with
blender-horseshoes\footnote{We refer to~\cite{BD12} to the precise
definition of $cu$ and $cs$-blender-horseshoe. Here it suffices to
understand that these objects are a certain class of horseshoes in
dimension $d\geq 3$.} which are more abundant objets than
hyperbolic
non-trivial attractors. %and folding manifold.}.
It is unknown if $C^1$-Newhouse domains can be obtained
arbitrarily {close to homoclinic tangency in dimension $d\geq 3$.}
It is also unknown whether $C^1$-Newhouse domains exist for
surface dynamics. Moreira's result~\cite{Moreira2011} provides a
strong evidence suggesting that there are no Newhouse domains in
the $C^1$-topology for surface dynamics.

\subsubsection{Newhouse domains of tangencies with {non-real complex leading multipliers}} First, we will consider Newhouse domains of
diffeomorphisms of dimension $m\geq 3$ with homoclinic tangencies
associated with periodic saddles whose multipliers $\lambda_1,
\lambda_2,\dots,\lambda_{m-1},\gamma$ satisfy that $\lambda_1$ and
$\lambda_2$ are non-real complex conjugate, that is,
$\lambda_{1,2}=\lambda e^{\pm i\varphi}$ with
$\lambda\in\mathbb{R}$, $\varphi\not=0,\pi$, and
\begin{equation}
\label{eq1}
   |\lambda_j|<|\lambda|<1<|\gamma| \ \ \ \text{with \ \ \ $|\lambda^2\gamma|<1<|\lambda\gamma|$} \quad \text{for $j\not= 1,2$.}
\end{equation}
Recall that a saddle is  \emph{sectional dissipative} if  the
product of any pair of multipliers is less than one in absolute
value. This implies that the unstable index (dimension of the
unstable manifold) needs to be one. A periodic point
satisfying~\eqref{eq1} has also unstable index one, it is
\emph{dissipative} (product of all multipliers is less than one)
but it is not sectional dissipative. {On the other hand, the \emph{leading multipliers} of a periodic point $P$ of period $n$ of a diffeomorphism $f$ are the nearest eigenvalues of $Df^n(P)$ to the unit circle. Thus, roughly speaking, we are considering Newhouse domains of homoclinic tangencies associated with dissipative but non-sectionally dissipative periodic points with non-real complex leading multipliers. Note that when we say~"non-real complex multipliers" we are asking that \emph{some} of the multipliers are not real (but not necessarily all).}

 As we have mentioned, this kind of Newhouse domains can
be obtained for $r\geq 2$ arbitrarily $C^r$-close to
diffeomorphism with a homoclinic tangency associated with a
periodic point satisfying~\eqref{eq1}. For $r=1$ (actually for any
$r\geq 1$), {let us consider} a
$C^r$-diffeomorphism $f$ with the following~properties:
%\begin{enumerate}[label=(H\arabic*)]
%\item \label{H1} $f$ has a equidimensional cycle between saddle periodic points
%$P$ and $Q$ where the multipliers of $Q$ satisfying~\eqref{eq1}.
%That is, both points have the same index (dimension of the stable
%bundle) and its invariant manifolds meet cyclically.
%\item \label{H2} $P$ belongs to a $cs$-blender-horseshoe $\Gamma$
%and $W^s(Q)$ and $W^u(P)$ has non-transverse intersection.
%\end{enumerate}
\begin{enumerate}[label=(H\arabic*)]
\item \label{H1} $f$ has a \emph{non-transverse equidimensional cycle} associated with hyperbolic periodic points  $P$ and $Q$.
That is, both $P$ and $Q$ have the same
unstable index, its stable and unstable invariant manifolds
meet {transversely and} cyclically  and also have at least one topologically non-transverse
intersection;
%\footnote{Since one of
%the invariant manifold is of codimension one, by topologically
%non-transverse here we understand that locally around the tangency
%point the one-dimensional invariant manifold belies in one of the
%pieces tha}.
\item \label{H3} $Q$ has multipliers satisfying~\eqref{eq1};
\item \label{H2} $P$ is homoclinically related to a $cs$-blender-horseshoe $\Gamma$.
\end{enumerate}

%Recall that two saddles of a diffeomorphism $f$ are homoclinically related if the invariant
%manifolds of their orbits intersect transversely and cyclically. To be homoclinic
%related defines an equivalence relation on the set of saddles of $f$. Two saddles
%that are homoclinically related have the same  unstable index. A saddle is homoclinically related to the blender
%if it is homoclinically related to some saddle of the blender.

{
It is not difficult to see that a $C^r$-Newhouse domain of homoclinic tangencies associated with periodic points satisfying~\eqref{eq1} can be also obtained arbitrarily $C^r$-close to $f$ under the assumptions~\ref{H1}-\ref{H3}-\ref{H2}. See Proposition~\ref{lem-new}.}
%This is immediately followed by using \cite[Sec.~4.3]{BD12} to get
%first a $C^1$-robust equidimensional tangency associated with the
%continuation of $Q$ and $\Gamma$ arbitrarily $C^r$-close to $f$.
%To be more precise, one gets an open set $\mathcal{N}$ of
%$C^r$-diffeomorphisms where $f\in \overline{\mathcal{N}}$ and
%every $g\in \mathcal{N}$ has a tangency between $W^s(Q_g)$ and
%$W^u(\Gamma_g)$ where $Q_g$ and $\Gamma_g$ are the continuations
%of $Q$ and $\Gamma$ for~$g$. After that, using the Inclination
%Lemma, any $g\in \mathcal{N}$ can approximate by a homoclinic
%tangency associated with the continuation~$Q$. Consequently,
%$\mathcal{N}$ is a $C^r$-Newhouse domain ($r\geq 1$) arbitrarily
%close to $f$.
As a corollary, we will obtain a similar result as
in~\cite{BD99} on the approximation of Newhouse domains associated
with homoclinic tangencies to saddle periodic points
satisfying~\eqref{eq1} from {heterodimensional cycles in dimension $d=3$ and non-transverse equidimensional cycle in dimension $d\geq 3$ with non-real complex leading multipliers.}
Moreover, we will show that any diffeomorphism in a
Newhouse domain of this type is contained in the closure of
diffeomorphisms having historic contractive non-trivial wandering
domains.

\begin{mainthm} \label{thm1} Let $\mathcal{N}$ be a $C^r$-Newhouse
domain ($r\geq 1$) of homoclinic tangencies associated with
periodic points satisfying~\eqref{eq1}. Then there is a dense set
$\mathcal{D}$ of $\mathcal{N}$ such~that, every $f\in \mathcal{D}$
has a non-trivial historic contractive wandering domain. Moreover,
the set $\mathcal{N}$ can be obtained arbitrarily
\begin{itemize}[leftmargin=0.5cm,itemsep=0.1cm]
\item $C^r$-close to a diffeomorphism having a
non-transverse equidimensional cycle  satisfying \ref{H1}-\ref{H3}
for $r\geq 2$ and satisfying~\ref{H1}-\ref{H3}-\ref{H2} for $r=1$.
\item $C^1$-close to a three-dimensional diffeomorphism having a heterodimensional
cycle associated  with a pair of hyperbolic periodic saddles with non-real
complex multipliers and where the multipliers of some of these
saddles satisfy~\eqref{eq1}.
\end{itemize}
\end{mainthm}

The idea behind the proof of the first part of the above theorem
is a reduction of the homoclinic tangency to a two-dimensional
smooth normally-hyperbolic attracting invariant manifold where the
restricted dynamics has a dissipative saddle. After that, we apply
the result in~\cite{KS17}. This strategy is not new and  was
successfully applied to find other types of complex dynamics
in~\cite{Ro95} and~\cite{KS06}. However, we cannot apply Romero's
result~\cite[Thm.~C]{Ro95}, even in the three-dimensional case
(see also~\cite[Thm.~4, Rem.~1, Sec.~4.1]{KS06} or
\cite[Lem.~2]{GTS93}) because of the following difficulty. Recall
first Romero's result in the three-dimensional case. Let  $f$ be a
$C^r$-diffeomorphism ($r\geq 2$) having a homoclinic tangency
associated with a periodic point $P$ with real multipliers
$\nu,\lambda,\gamma$ such that
$$
    |\nu|<|\lambda|<1<|\gamma| \quad \text{and} \quad J(P)\eqdef |\lambda
    \gamma|>1.
$$
The case where $P$ has complex multipliers is reduced to the above
case using~\cite[Sec.~5]{PV94}. According to~\cite[Thm.~C]{Ro95},
arbitrarily $C^r$-close to $f$ there exists a diffeomorphism $g$
which has a two-dimensional normally hyperbolic attracting smooth
invariant manifold $S$ such that the two-dimensional restriction
$g|_{S}$ has a homoclinic tangency associated with a
periodic point $Q$ with $J(Q)>1$. Since $Q$ is not a dissipative
periodic point, we cannot apply~\cite{KS17} to $g|_S$. To work around
this problem, we use the rescaling theory in~\cite{GST08} working
directly with the complex multipliers instead of reducing the
problem to the case of real leading multipliers as
in~\cite{PV94,Ro95}. %\\[-0.7cm]

\subsubsection{Historic wandering domains from Tatjer homoclinic tangencies}

In dimension three, we can also obtain wandering domains from
another type of Newhouse domains. Namely, we will consider
Newhouse domains associated with Tatjer homoclinic tangencies.
To introduce these tangencies, we need some
preliminaries.

Let $P$ be a hyperbolic saddle fixed point of a three-dimensional
diffeomorphism $f$. For simplicity of the exposition, we have
chosen a fixed point, but all terminologies and concepts are
valid if $P$ is a periodic point. Suppose that $Df(P)$ has real
eigenvalues $\lambda_{s}$, $\lambda_{cu}$ and $\lambda_{uu}$
satisfying  \[
   |\lambda_{s}| < 1< |\lambda_{cu}|<|\lambda_{uu}|.
\] Thus the tangent space at $P$ has a dominated splitting of the
form $E^{s}\oplus E^{cu} \oplus E^{uu}$ given by the corresponding
eigenspaces. The unstable manifold $W^{u}(P)$ is tangent at $P$ to
the bundle $E^u=E^{cu}\oplus E^{uu}$. On the other hand, according
to~\cite{HPS77}, the extremal bundle $E^{uu}$ can be also
integrated providing a one-dimensional manifold $W^{uu}(P)$ called
strong unstable manifold. Moreover, this bundle can be uniquely
extended to $W^u(P)$ providing a foliation $\mathcal{F}^{uu}(P)$ of
this manifold by one-dimensional leaves $\ell^{uu}(Y)$ containing
$Y\in W^u(P)$. We assume additionally that the center-stable
bundle $E^{s}\oplus E^{cu}$ is also extended and integrated along
the stable manifold $W^s(P)$ of $P$. Although the extended
center-stable bundle is not unique, any center-stable manifold
contains $W^s(P)$ and any two of these  manifolds are tangent to
each other at every point of $W^u(P)$. Finally,
%\begin{defi}
a three-dimensional diffeomorphism as above has a \emph{Tatjer
homoclinic tangency}  associated with $P$ (which corresponds to
the type~I in~\cite{Ta01})~if
\begin{enumerate}[leftmargin=1.2cm,rightmargin=0.8cm,label=(T\arabic*)]
\item \label{Tatjer1} $W^s(P)$ and $W^u(P)$ have a quadratic tangency at $Y$ which
does not belong to the strong unstable manifold $W^{uu}(P)$ of
$P$,
\item \label{Tatjer2} $W^s(P)$ is tangent to the leaf $\ell^{uu}(Y)$ of
$\mathcal{F}^{uu}(P)$ at $Y$,
\item \label{Tatjer3} $W^u(P)$ is transverse to any center-stable manifold at
$Y$.
\end{enumerate}
%\end{defi}
If $P$ has stable index  equals two, the above definition
%of Tatjer homoclinic tangency
applies to $f^{-1}$.

%We must also notify that originally Tatjer in~\cite{Ta01} includes
%the extra assumption of $C^1$-linearazing coordinates around $P$.
%Later in~\cite{gonchenko2007bifurcations} the results
%in~\cite{Ta01} were generalized without this assumption.
Similarly to the results obtained in~\cite{GST93b,GST08}
%for the
%unfolding of homoclinic tangencies associated with periodic points
%with complex multipliers as in~\eqref{eq1},
strange attractors, normally hyperbolic attracting smooth
invariant circles and hyperbolic sinks are also obtained by
unfolding a Tajter homoclinic tangency under the following extra
assumptions~\cite{Ta01,GGT07}. The first extra
assumption is the \emph{dissipativeness}:  the homoclinic tangency
is associated with a saddle periodic point $P$ whose multipliers
are $\lambda_{s}$, $\lambda_{c}$ and $\lambda_{u}$ with
\[
|\lambda_s|<1<|\lambda_u|, \quad
|\lambda_{s}|<|\lambda_c|<|\lambda_u| \quad \text{and} \quad
   |\lambda_{s} \lambda_{c} \lambda_{u}|<1  \ \
   \text{(dissipativeness).}
\]
Recall a periodic point is said to be sectional dissipative when
the absolute value of the product of any pair of multipliers is less than one.
Conversely, the second extra assumption required is the
\emph{non-sectional dissipativeness} of $P$:  either
\begin{align*} %\label{autovalorescaseA}
  &\text{(Case A)}  \qquad \qquad
   |\lambda_c|<1, \qquad|\lambda_{c}\lambda_{u}|>1  %\
 %  \text{or}
   \\
 %\label{autovalorescaseB}
  &\text{(Case B)}  \qquad \qquad \ |\lambda_{c}|>1.
\end{align*}
For short we will say that $P$ is \emph{dissipative but
non-sectional dissipative periodic point} when both above extra
assumptions are satisfied.

On the other hand, observe that the conditions~\ref{Tatjer1} and
\ref{Tatjer3} are generic. This means that for an arbitrarily small
perturbation one can always assume that a homoclinic tangency
under the assumption~\ref{Tatjer2} is, in fact, a Tatjer tangency
(of type~I). Although~\ref{Tatjer1} is a codimension one
condition, we must observe that the required tangency
in~\ref{Tatjer2} is a condition of codimension
\[
3-\dim [T_YW^u(P)+T_Y\ell^{uu}(Y)]=2.
\]
For more details about tangencies of large codimension see
also~\cite{BR17,BR21,BP20}.

\begin{mainthm} \label{thm2}
Let $f$ be a three-dimensional $C^r$-diffeomorphism ($r\geq 2$)
with a Tatjer homoclinic tangency associated with a dissipative
but non-sectional dissipative periodic point. Then,
$C^r$-arbitrarily close to $f$, there are a $C^r$-Newhouse domain
$\mathcal{N}$ (associated with sectional dissipative periodic
points) and a dense subset $\mathcal{D}$ of $\mathcal{N}$ such
that, every $f\in \mathcal{D}$ has a non-trivial historic
\mbox{contractive wandering domain.}
% there exist
%arbitrarily $C^r$-close to $f$ a $C^r$-Newhouse domain
%$\mathcal{N}$ and a dense subset $\mathcal{D}$ of $\mathcal{N}$
%such that, every $f\in \mathcal{D}$ has a historic contracting
%non-trivial wandering domain.
%For any $r\geq 1$, there exists  a $C^r$-Newhouse domain
%$\mathcal{N}$ associated with Tatjer homoclinic tangencies
%satisfying the dissipativeness condition and either Case~A or
%Case~B condition. Moreover, there are a dense subset $\mathcal{D}$
%of $\mathcal{N}$ such that, every $f\in \mathcal{D}$ has a
%non-trivial historic contractive wandering domain. And,
%furthermore, $C^r$-Newhouse domains associated with Tatjer
%homoclinic tangencies can be obtained arbitrarily $C^1$-close to a
%3-dimensional diffeomorphisms having a heterodimensional cycle
%associated with a pair of of hyperbolic periodic saddles with
%complex multipliers.
\end{mainthm}

The following result completes the case $r=1$. But to
achieve this, we need to introduce a special type of Newhouse
domains. Namely, we will deal with a Newhouse domain $\mathcal{N}$
satisfying  that there exists a dense set $\mathcal{D}$ of
$\mathcal{N}$ such that a map in $\mathcal{D}$ displays a Tatjer
homoclinic tangency associated with dissipative but non-sectional
dissipative periodic point. Observe that because of the extra
degeneration on the codimension of the homoclinic tangency, the
existence of a Newhouse domain associated with Tatjer tangencies
 is a non-trivial problem. However, as we will explain in~\S\ref{sec:proofC},
 following essentially
the strategy of~\cite{KNS17}, we construct such Newhouse domains
near certain non-transverse equidimensional cycles. Finally, we
will show that this open class of diffeomorphisms also is in the
closure of maps having non-trivial historic contractive wandering
domains.

\begin{mainthm} \label{thm3}
Let $\mathcal{N}$ be a $C^r$-Newhouse domain of Tatjer homoclinic
tangencies associated with  dissipative but non-sectional
dissipative periodic points with $r\geq 1$. Then, there is a dense
set $\mathcal{D}$ of $\mathcal{N}$ such~that, every $f\in
\mathcal{D}$ has a non-trivial historic contractive wandering
domain.
%Moreover, these type of $C^r$-Newhouse domains can be obtained
%$C^1$-arbitrarily close to a $C^r$-diffeomorphism having
%heterodimensional cycle associated with periodic points with
%complex multipliers where at least one of them is non-sectional
%dissipative.

Moreover, this type of $C^r$-Newhouse domains can be obtained
arbitrarily  $C^1$-close to a $C^r$-diffeomorphism having
heterodimensional cycle associated with periodic points with non-real
complex multipliers such that at least one of them is dissipative
but non-sectional dissipative.
\end{mainthm}

%Let $f$ be a $C^r$-diffeormorphisms of a 3-dimensional manifold
%with a heterodimensional cycle associated with periodic points $P$
%and $Q$ with complex multipliers. Assuming that the multiplies of
%one of the periodic point is dissipative but not sectional
%dissipative (i.e., if $\lambda_1,\lambda_2,\lambda_3$ are the
%multipliers then either
%$$
%   |\lambda| < 1 <|\gamma|  \ \ \ \text{and} \ \ \
%   |\lambda^2 \gamma|<1<|\lambda \gamma| \ \ \ \text{where
%   $\lambda_{1,2}=\lambda e^{i\varphi}$ and $\gamma=\lambda_3$
%   with $\varphi\not= 0,\pi$}
%$$
%or
%$$
%  |\lambda| < 1 <|\gamma|  \ \ \ \text{and} \ \ \
%   |\lambda \gamma^2|<1 \ \ \ \text{where
%   $\lambda=\lambda_1$ and  $\lambda_{2,3}=\gamma e^{i\varphi}$
%   with $\varphi\not= 0,\pi$}.
%$$
%%\section{Proof of the theorems}

{As we explain in Section~\ref{sec:proofC}, the proof of the main result in~\cite{KNS17} has a gap. The additional assumption on the multipliers that appears in the second part of the above theorem solves this gap. Actually, the conclusion of Theorem~\ref{thm3} is stronger than~\cite[Thm.~1.1]{KNS17}  since the non-trivial contracting wandering domains  obtained are also~\emph{historic}.}

%\vspace{0.1cm}
\subsection{Attracting circles, strange attractors, sinks and non-trivial homoclinic classes}
%Mentioned that,
From \cite{GST08}, it also follows the coexistence
of infinitely many normally hyperbolic attracting invariant smooth
circles (and sinks) for a residual subset of diffeomorphisms in a
$C^r$-Newhouse domain associated with saddle periodic points
satisfying~\eqref{eq1}. This result is only proved in~\cite{GST08}
in the case $r\geq 2$ but the case $r=1$ also holds since
$C^\infty$-diffeomorphisms with homoclinic tangencies are
$C^1$-dense in a $C^1$-Newhouse domain and the attracting circles
(and sinks) are $C^1$-robust (they are normally hyperbolic).
From~\cite{Ta01,GGT07}, the same results are obtained for
$C^r$-Newhouse domains ($r\geq 1$) of Tatjer homoclinic tangencies
associated with dissipative but non-sectional dissipative periodic~points.

On the other hand, attracting compact invariant sets having a
dense orbit with at least one positive Lyapunov exponent obtained
from H\'enon-like maps, the so-called \emph{H\'enon-like strange
attractors}, are non-hyperbolic. This lack of hyperbolicity
prevents stability under perturbations, and thus, the classical
arguments (see~\cite{PT93}) to provide coexistence of infinitely
many of such attractors do not work. This difficulty was overcome
by Colli~\cite{colli1998infinitely} and Leal~\cite{leal2008high}.
From these papers, it follows that, in a $C^r$-Newhouse domain ($r\geq 1$)
associated with homoclinic tangencies of sectional dissipative
periodic points, there exists a dense set of diffeomorphisms
exhibiting the coexistence of infinitely many non-hyperbolic
strange attractors (see~\S\ref{sec:proofD}). Once again, this
result can be translated to the $C^r$-Newhouse domains considering
in Theorems~\ref{thm1} and \ref{thm3}. This is because the main
tool behind the proof of  these theorems is a reduction of the
dynamics to a two-dimensional attracting smooth invariant manifold
where the restriction of the diffeomorphism has a homoclinic
tangency associated with a dissipative periodic point. Then, one
can apply~\cite{KNS17} as well as~\cite{colli1998infinitely}.

%\vspace{0.3cm}

Notice that H\'enon-like strange attractors are, in fact,
non-trivial attracting homoclinic classes. Recall that a
\emph{homoclinic class} is the closure of the transverse
intersections of the invariant manifolds (stable and unstable
ones) of the hyperbolic periodic orbit. By attracting we
understand that there exists an open neighborhood $V$ of the
homoclinic class such that the forward image of the closure of $V$
is strictly inside of $V$. And, we say that the homoclinic class
is \emph{non-trivial} if it is not reduced to a sink or repeller.
Although, as mentioned, H\'enon-like strange attractors are not
stable under perturbation, non-trivial attracting homoclinic
classes are $C^1$-robust. This observation allows us to get the
following result:

\begin{mainthm} \label{thm4}
Let $\mathcal{N}$ be a $C^r$-Newhouse domain ($r \geq 1$) of one
of the following types:
\begin{enumerate}[label=(\roman*)]
\item \label{1} homoclinic tangencies associated with sectional dissipative periodic
points,
\item \label{2} homoclinic tangencies associated with periodic points
satisfying~\eqref{eq1},
\item \label{3}  Tatjer
homoclinic tangencies associated with dissipative but \\
non-sectional  dissipative periodic points. %\\[-0.7cm]
\end{enumerate}
Then, there is a residual set $\mathcal{R}$ of $\mathcal{N}$
such~that, every $f\in \mathcal{R}$ exhibits the coexistence of
infinitely many (independent) non-trivial attracting homoclinic
classes.
\end{mainthm}

In the topology $C^1$, it has been known for sometime~\cite{BD03}
that the coexistence of infinitely many (pairwise disjoint)
non-trivial attracting homoclinic classes are locally generic. In
fact, recently, it has been also proved~\cite[Thm.~5]{BCF18} that
for $C^1$-generic diffeomorphisms these homoclinic classes could
be taken with entropy uniformly large. To conclude, we want to
remark that, also in the $C^1$-topology,  %as a consequence of the previous theorems
we obtain the following:

\begin{maincor} Arbitrarily $C^1$-close to a three-dimensional
$C^r$-diffeomorphism $f$ having a heterodimensional cycle associated
with periodic points with non-real complex multipliers where at least one
of them is dissipative but non-sectional dissipative, there exists
a locally residual set of diffeomorphisms exhibiting the
coexistence of infinitely many (independent) non-trivial
attracting homoclinic classes.
\end{maincor}

{The proof of this result is just applied Theorem~\ref{thm4} to the $C^r$-Newhouse domain of Tatjer
homoclinic tangencies associated with dissipative but
non-sectional  dissipative periodic points obtained from Theorem~\ref{thm3}  arbitrarily $C^1$-close to $f$.}

%\vspace{-0.5cm}
\section{Proof of the theorems}

{ Before proving the main theorems, let us provide the formal definition of the main two dynamical configurations in this paper:

\begin{defi}  \label{def-cycles}
  A diffeomorphism $f$ has a
  \begin{enumerate}[leftmargin=0.75cm]
    \item \emph{homoclinic tangency} associated with a transitive hyperbolic set $\Lambda$ if  there is a pair of points $x,y\in \Lambda$ such that $W^s(x)$ and $W^u(y)$ has a non-transverse intersection. The tangency is said to be \emph{$C^1$-robust} if there is a $C^1$-neighborhood $\mathcal{U}$ of $f$ such that any $g\in \mathcal{U}$ has a homoclinic tangency associated with the continuation $\Lambda_g$ of $\Lambda$ for $g$.
    \item \emph{heterodimensional cycle} associated with transitive hyperbolic sets $\Lambda$ and $\Sigma$ if these sets have different unstable indices and their invariant manifolds meet cyclically, that is,  $$W^s(\Lambda)\cap W^u(\Sigma)\not=\emptyset \quad \text{and} \quad W^u(\Lambda)\cap W^s(\Sigma)\neq\emptyset.$$
          The heterodimensional cycle is said to be \emph{$C^1$-robust} if there is a $C^1$-neighborhood $\mathcal{U}$ of $f$ such that any $g\in \mathcal{U}$ has a heterodimensional cycle associated with the continuations $\Lambda_g$ and $\Sigma_g$ of $\Lambda$ and $\Sigma$ for $g$ respectively.
  \end{enumerate}
\end{defi}
}

\subsection{Proof of Theorem~\ref{thm1}} \label{sec:proofA}
Let us assume that $\mathcal{N}$ is a $C^r$-Newhouse domain of
homoclinic tangencies associated with periodic points
satisfying~\eqref{eq1} with $r\geq 1$. Recall, according
to~\cite{GST08}, a homoclinic tangency is said to be \emph{simple}
if the tangency is quadratic, of codimension one and, in the case
that the dimension {$m\geq3$}, any extended unstable manifold is
transverse to the leaf of the strong stable foliation which passes
through the tangency point. Thus, since these properties are
generic, by an arbitrarily small $C^r$-perturbation with $r\geq
1$, we obtain that maps $f$ with a simple homoclinic tangency
associated with a periodic point $Q$ satisfying~\eqref{eq1} can be
obtained densely in the $C^r$-Newhouse domain $\mathcal{N}$.
Moreover, we can assume that $f$ is in fact $C^k$ with $k> r$.
%Proof: $f$ can be aproximated by $C^\infty$-diffeomorphisms.
%Then we can take take a family $g_\varepsilon$ of $C^\infty$ maps
%with $g_0$ arbitrarily  close $f$. This family looks like
%a generic  unfolding of $f$. Thus, parameters having
%new homoclinic tangencies associated with periodic points
%homoclinically related with the continuation of $Q$ appears (see for instance~\cite[sec.~5]{PV94}).
%Then, applying the inclination Lemma,
%we can find a homoclinic tangency of a $C^\infty$-diffeomorphims
%associated with the continuation of $Q$.

We need to consider a two-parameters unfolding $f_\varepsilon$ of
$f=f_0$ with $\varepsilon=(\mu,\varphi)$ where $\mu$ is the
parameter that controls the splitting of the tangency and
$\varphi$ is the value for which the argument of the complex
multiplier of $Q$ is perturbed. As usual, $T_0=T_0(\varepsilon)$
denotes the local map. In this case, this map corresponds to
$f^q_{\varepsilon}$, where $q$ is the period of $Q$ and it is
defined on a neighborhood $W$ of $Q$. By $T_1=T_1(\varepsilon)$ we
denote the map $f_\varepsilon^{n_0}$ from a neighborhood $\Pi^-$
of a tangent point $Y^-\in W^u_{loc}(Q,f_0)\cap W$ of $f_0$ to a
neighborhood $\Pi^+$ of $Y^+=f_0^{n_0}(Y^-)\in
W^s_{loc}(Q,f_0)\cap W$. Then, for $n$ large enough, one defines
the first return map $T_n=T_1\circ T_0^n$ on a subset
$\sigma_n=T_0^{-n}(\Pi^-)\cap \Pi^+$ of $\Pi^+$ where $\sigma_n\to
W^s_{loc}(Q)$ as $n\to\infty$.  According to~\cite[Lemma~1
and~3]{GST08} we have the following result:

\begin{lem}\label{lema-GHM}
There exists a sequence of open set $\Delta_n$ of parameters
converging to $\varepsilon=0$ such that for these values the map
$T_n$ has a two-dimensional normally hyperbolic attracting
invariant $C^k$-manifold $\mathcal{M}_n$ in $\sigma_n$ which,
after a $C^k$-smooth transformation of coordinates on $\sigma_n$,
the restriction of the map is given by
\begin{equation} \label{eq-GHM}
   \bar{x}=y, \qquad \bar{y}=M-Bx-y^2-R_n(xy+ o(1)).
\end{equation}
 The rescaled parameters $M$, $B$ and $R_n$ are functions of
$\varepsilon \in \Delta_n$ such that $R_n$ converges to zero as
$n\to \infty$ and $M$ and $B$ run over asymptotically large
regions which, as $n\to \infty$, cover all finite values. Namely,
\begin{align*}
M \sim \gamma^{2n}(\mu + O(\gamma^{-n}+\lambda^n)), \quad B \sim
(\lambda \gamma)^{n}\cos(n\varphi+o(1)) \quad \text{and} \quad R_n
\sim \frac{2J_1}{B}(\lambda^2\gamma)^n
\end{align*}
where $J_1\not=0$ is the Jacobian of the global map $T_1$
calculated at the homoclinic point $Y^-$ for $\varepsilon=0$. The
$o(1)$-terms tend to zero as $n\to \infty$ along with all the
derivatives up to the order $k$ with respect to the coordinates
and up to the order $k-2$ with respect to the rescaled parameters
$M$ and $B$. %Moreover, the limit family is the H\'enon map.
\end{lem}

%\vspace{0.3cm}

\begin{figure}
\labellist \small\hair 2pt \pinlabel $\mathrm{SN}_n$ at 400 390
\pinlabel $\mathrm{H}_n$ at 176 650 \pinlabel $\mathrm{T}_n^+$ at
362 590 \pinlabel $\mathrm{BT}_n$ at 250 440 \pinlabel
{$\mathrm{T}^-_n$} at 298 575
\endlabellist
  \includegraphics[scale=0.8]{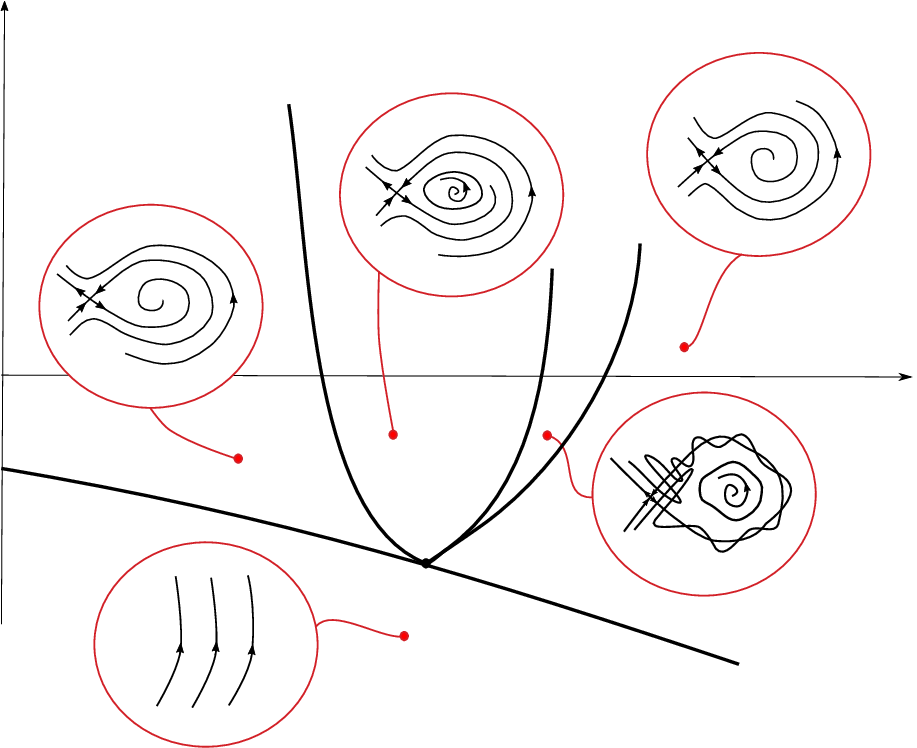}
\caption{Bifurcation diagram near the Bogdanov-Takens point
$\mathrm{BT}_n$ in the cases $R_n > 0$. The case $R_n<0$ is
similarly changing the stability of the periodic points. The curves
$\mathrm{SN}_n\setminus \{\mathrm{BT}_n\}$ and $\mathrm{H}_n$
correspond to saddle-node and Hopf bifurcations. The curves
$\mathrm{T}^-_n$ and $\mathrm{T}^+_n$ are curves of homoclinic
tangencies associated with a dissipative fixed point.}
\label{fig-BT}
\end{figure}
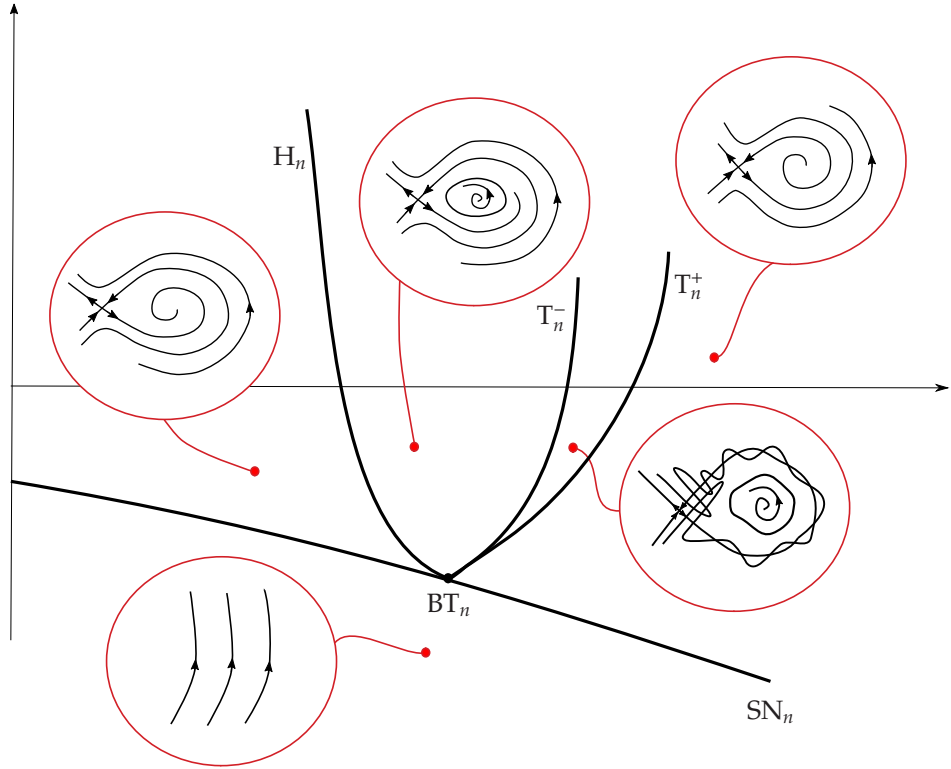

The dynamics of the generalized H\'enon map
\begin{equation} \label{eq-GHM-simplified}
   \bar{x}=y, \qquad \bar{y}=M-Bx-y^2-R_nxy
\end{equation}
was studied in~\cite{GG00,GG04,GKM05} (see also \cite{GGT07}). For
small $R_n$, the map~\eqref{eq-GHM-simplified} has, on the
parameter plane $(M,B)$, a bifurcation point
\begin{equation} \label{eq:pontos}
\begin{aligned}
\mathrm{BT}_n \ &: \quad  M = \frac{-1 -R_n}{(1 + R_n/2)^2}, \quad
B=1+\frac{R_n}{1+R_n/2}.
\end{aligned}
\end{equation}
At this point, \eqref{eq-GHM-simplified} has a fixed point with a
pair of eigenvalues {equal} to $+1$.  As it was shown
in~\cite{GKM05} (see also~\cite{GGT07,Ta01}), the Generalized
H\'enon family unfolds generically a~\emph{Bogdanov-Takens
bifurcation} at $\mathrm{BT}_n$. Figure~\ref{fig-BT} is showed
the local picture of this bifurcation
(c.f.~\cite{broer1996invariant}).

%\vspace{0.2cm}

Although the coefficient $R_n$ in~\eqref{eq-GHM} depends on $B$,
note that the range of values it takes is negligible when $B$ is
limited and $n$ is large enough. Thus, the bifurcation diagram
of~\eqref{eq-GHM} can be studied from the results described above
for~\eqref{eq-GHM-simplified} assuming $R_n=o(1)$ independent of
$B$. Thus, for any $n$ large enough, there are values of the
parameter $\varepsilon \in \Delta_n$ such that the parameters
$M=M(\varepsilon)$ and $B=B(\varepsilon)$ of
$T_n=T_n(\varepsilon)$ belong to the curves $T_n^\pm$ in
Figure~\ref{fig-BT}.  Thus, $T_n$ has at {these} parameters values a homoclinic tangency associated with a dissipative saddle fixed point. In other words, we can find a sequence $(\varepsilon_n)_n$ of parameters $\varepsilon_n \in \Delta_n$ with $\varepsilon_n\to 0$ such that $g_n=f_{\varepsilon_n}$ (which approaches $f_0=f$) has a normally hyperbolic attracting two-dimensional invariant
manifold $\mathcal{M}_n$ for some iterated $m=m(n)$ where the
restriction of $g_n^m$ to this manifold has a homoclinic tangency
associated with a dissipative periodic point. Now,
applying~\cite{KS17} to the restriction
$g_n^{m}|_{\mathcal{M}_n}$, we obtain a map arbitrarily
$C^k$-close to $g_n$ with a non-trivial historic contractive
wandering domain. In particular, we obtain that the set of maps
with wandering domains are $C^r$-dense in $\mathcal{N}$. This
completes the first part of Theorem~\ref{thm1}.

To prove the second part, consider firstly the case of a
non-transverse equidimensional cycle (assumption~\ref{H1}). As it
is well-known using the Inclination Lemma~\cite[Lemma~7.1]{junior1982geometric}, a $C^r$-diffeomorphism
having a non-transverse equidimensional cycle can be
$C^r$-approximated by diffeomorphisms exhibiting homoclinic
tangencies. Moreover, from assumption~\ref{H3}, these homoclinic
tangencies can be obtained associated with periodic points
satisfying~\eqref{eq1}. Hence, as it was mentioned in the
introduction, according to~\cite{GTS93} one can obtain a
$C^r$-Newhouse domain as desired for $r\geq 2$ arbitrarily close
to a diffeomorphism exhibiting homoclinic tangencies associated
with periodic points satisfying~\eqref{eq1}. {Next proposition proves the case $r=1$. First, we need some definitions.}

{Recall that two periodic points of $f$
 are {\emph{homoclinically related}} if the invariant manifolds of their orbits intersect transversely and cyclically.
 To be homoclinic related defines an
equivalence relation on the set of periodic points of $f$. Two saddles that are homoclinically related have the same
unstable index. We can generalize the above notion for a pair of transitive hyperbolic sets $\Lambda$ and $\Sigma$ of $f$ by saying that they are homoclinically related if there is a pair of periodic points $P \in \Lambda$ and $Q\in \Sigma$ homoclinically related. Note that in a transitive hyperbolic set all  pairs of periodic points are homoclinically related. Thus, the above periodic points $P$ and $Q$  are homoclinically related with any other periodic point $R \in \Lambda\cup \Sigma$.}

{We say that   diffeomorphism $f$ has a \emph{non-transverse equidimensional cycle} associated with transitive hyperbolic sets $\Lambda$ and $\Sigma$
if $\Lambda$ and $\Gamma$ are homoclinically related and have a \emph{heterodimensional tangency}, i.e., there are $x\in \Lambda$ and $y\in \Sigma$ such that $W^s(x)$ and
$W^u(y)$ has tangency. Moreover, a non-transverse equidimensional cycle as above is said to be \emph{$C^1$-robust} if there exists a $C^1$-neighborhood $\mathcal{U}$ of $f$  such that every $g\in \mathcal{U}$ has as non-transverse equidimensional cycle associated with the continuations $\Gamma_g$ and $\Sigma_g$ of $\Lambda$ and $\Sigma$ respectively. }

{
\begin{prop} \label{lem-new}
Let $f$ be a $C^r$-diffeomorphism ($r\geq 1$) under the assumptions~\ref{H1} and \ref{H2}. Then $f$ can be
$C^r$-approximated by diffeormorphisms exhibiting a $C^1$-robust
non-transverse equidimensional cycle associated with the
continuations of the periodic point $Q$ and the
$cs$-blender-horseshoe $\Gamma$.

Moreover, if in addition $f$ satisfies also~\ref{H3}, then  a $C^r$-Newhouse domain of homoclinic tangencies associated with
%dissipative but non-sectionally dissipative
periodic points
%with complex leading multipliers (i.e.,
satisfying~\eqref{eq1}  %)
is also obtained arbitrarily $C^r$-close to $f$.
\end{prop}

\begin{proof}
By assumption~\ref{H1}, $P$ and $Q$ are homoclinically related and
by~\ref{H2}, $P$ is homoclinically related to a $cs$-blender-horseshoe $\Gamma$.
 Thus $P$ and $Q$ are homoclinically related to any saddle of $\Gamma$ and, in particular, to the so-called reference saddle of the blender-horseshoe (see~\cite[Def.~3.9]{BD12}).  Now, by assumption~\ref{H1}, $P$ and $Q$ have a heterodimensional tangency. Using the Inclination Lemma~\cite[Lemma~7.1]{junior1982geometric}, one obtain easily  a diffeomorphism $g$ arbitrarily $C^r$-close to $f$ having a tangency between the reference saddle of $\Gamma_g$  and the stable manifold of $Q_g$. Here $\Gamma_g$ and $Q_g$ are the continuations of $\Gamma$ and $Q$ for $g$ respectively.  Then, according to~\cite[Lemma~4.10 and Corollary~4.11]{BD12}, we immediately conclude the first part of the proposition.

The second part of the proposition is just again a standard application of the Inclination
Lemma. Indeed, from the first part of the lemma, we get a diffeomorphism $g$ arbitrarily $C^r$-close to $f$ with a $C^1$-robust equidimensional cycle associated with  $Q_g$ and $\Gamma_g$. This provides a $C^r$-open set $\mathcal{N}_g$ where each diffeomorphism $h\in \mathcal{N}_g$  has a heterodimensional tangency between the invariant manifolds of $Q_h$ and a point $x\in \Gamma_h$. Since $x$ and $Q_h$ are also homoclinically related, using the Inclination Lemma, we obtain arbitrarily $C^r$-close to $h$ (in particular in $\mathcal{N}_g$) a new diffeomorphism  having a homoclinic tangency associated with the continuation of~$Q$. Consequently,
$\mathcal{N}_g$ is a $C^r$-Newhouse domain ($r\geq 1$) of homoclinic tangencies arbitrarily
$C^r$-close to $f$. Finally, since the continuation of $Q$ satisfies~\ref{H3},  then these homoclinic tangencies  are associated with periodic points satisfying~\eqref{eq1}. This completes the proof of the proposition.
\end{proof}
}

%On the other hand,
%for $r=1$, as we also discussed in the introduction,
%assumptions~\ref{H1}-\ref{H3}-\ref{H2} imply that $f$ can be
%approximated by diffeormorphisms displaying a $C^1$-robust
%non-transverse equidimensional cycle associated with the
%continuations of the periodic point $Q$ and the
%$cs$-blender-horseshoe $\Gamma$. In particular, arguing as above
%one can get a homoclinic tangency from an arbitrarily small
%perturbation of the cycle. This proves that $C^1$-Newhouse domains
%associated with periodic points satisfying~\eqref{eq1} can be
%obtained arbitrarily $C^1$-close to non-transverse equidimensional
%cycles under the assumptions~\ref{H1}-\ref{H3}-\ref{H2}.

Now we will prove that a $C^1$-Newhouse domain associated with
periodic points satisfying~\eqref{eq1} can be obtained arbitrarily
$C^1$-close to a three-dimensional diffeomorphism exhibiting a certain type of
heterodimensional cycle. Namely, a heterodimensional cycle associated  with a pair of hyperbolic
periodic points  with non-real complex multipliers and where some of them
satisfies~\eqref{eq1}. {To prove this, in view of Proposition~\ref{lem-new}, it suffices the following result:}

{
\begin{prop} \label{prop-H123}
Let $f$ be a three-dimensional $C^r$-diffeomorphism exhibiting a
heterodimensional cycle associated  with a pair of hyperbolic
periodic points $P$ and $Q$ with non-real complex multipliers.  Then $f$ can
be $C^1$-approximated by $C^r$-diffeormorphisms exhibiting a
non-transverse equidimensional cycle under the
assumptions~\ref{H1} and \ref{H2}.
Moreover, if in addiction $Q$ satisfies~\eqref{eq1}, then the non-transverse equidimensional cycle also verifies~\ref{H3}.
\end{prop}

According to~\cite[Prop.~2.1]{KNS17} we can $C^1$-approximate $f$
by diffeomorphisms exhibiting a non-transverse equidimensional
cycle under the assumption~\ref{H1}. Actually, the
first step in~\cite[Lemma~2.2]{KNS17} to prove this result was to
approach $f$ {in the $C^1$-topology} by a $C^r$-diffeomorphism $h$ having a heterodimensional cycle
associated with a pair of periodic points {$P'_h$ and $Q'_h$} with real multipliers
which are homoclinically related to {the continuations $P_h$ and $Q_h$ of} $P$ and $Q$ respectively
(c.f.~\cite[Thm.~2.1]{BD08}).  {We will show that a stronger result can be obtained. Namely, we will explain that, in addition, $P'_h$ is homoclinically related to a $cs$-blender-horseshoe. Thus~\ref{H2} holds. This will be obtained as an immediate application of the following general lemma for the so-called co-index one heterodimensional cycles of diffeomorphisms in any dimension. First, we need to introduce some notation.}

Denote by $\uindex(R)$ the unstable index of a hyperbolic periodic orbit $R$. Following~\cite[Def.~1.3]{BDK12}, a periodic point has \emph{real central multipliers} if \emph{all} the leading multipliers are real and with multiplicity one.

}

{
\begin{lem} \label{lem-medio} Let $f$ be a $C^r$-diffeomorphism having a  heterodimensional cycle associated  with a pair of hyperbolic
periodic points $P$ and $Q$ with $\uindex(Q) = \uindex(P)+1$ and where  at least one
of them has non-real complex leading multipliers. Then, arbitrarily $C^1$-close to $f$, there is a $C^r$-diffeomorphism $h$ having
\begin{enumerate} \item a $cs$-blender-horseshoe $\Gamma^{cs}_h$ and a $cu$-blender-horseshoe $\Gamma^{cu}_h$,
\item  periodic points $P'_h$ and $Q'_h$ with real central multipliers and homoclinically related to  the continuations $P_h$ and $Q_h$ of $P$ and $Q$ respectively,
\end{enumerate}
such that
\begin{enumerate}[resume]
\item $h$ has a heterodimensional cycle
associated with $P'_h$ and $Q'_h$,
\item $P'_h$ and $Q'_h$ are homoclinically related to $\Gamma^{cs}_h$ and $\Gamma^{cu}_h$ respectively.
\end{enumerate}
\end{lem}

\begin{proof}
%In order to fit the assumptions of the lemma with the referred literature, we must consider $f^{-1}$ and prove the lemma for such inverse. Thus, we have to show that arbitrarily $C^1$-close to $f^{-1}$, there is a $C^r$-diffeomorphism $h$ having $cu$-blender-horseshoe $\Gamma_h$ and a hete\-ro\-di\-men\-sional cycle associated with a pair of periodic points $P'_h$ and $Q'_h$ with real leading multipliers
%which are homoclinically related to the continuations $P_h$ and $Q_h$ of $P$ and $Q$ respectively and where $Q'_h$  is homoclinically related to $\Gamma_h$.

Using the stabilization theory
in~\cite[Thm.~1 and 2]{BDK12}, we find a $C^r$-diffeomorphism~$h$
arbitrarily $C^1$-close to $f$ having a $C^1$-robust
heterodimensional cycle between transitive hyperbolic sets $\Lambda_h$ and $\Sigma_h$
containing the continuation $P_h$ and $Q_h$ of $P$ and $Q$ respectively. Now, we have to go into the proof of this result to explain that, actually,  $h$ has a heterodimensional cycle associated with a pair of periodic points $P'_h$ and $Q'_h$ with real central multipliers belonging to $\Lambda_h$ and $\Sigma_h$ respectively. Moreover, $P'_h$ and $Q'_h$  belong to a $cs$-blender-horseshoe~$\Gamma^{cs}_h$ and $cu$-blender-horseshoe $\Gamma^{cu}_h$ respectively. This will prove the lemma.

{The key dynamical configuration is described in~\cite[Prop.~6.1]{BDK12} and depicted in \cite[fig.~6]{BDK12}. In this configuration, we have a strong homoclinic intersection associated with a partially hyperbolic saddle-node/flip $S_g$ whose strong stable and unstable manifolds also meet the invariant manifolds of $P_g$ and $Q_g$. Note that this configuration is achieved under the assumption that the initial heterodimensional cycle is \emph{non-twisted} (see~\cite[Fig.~1]{BDK12} for an intuitive explanation of these sort of heterodimensional cycles). However, this is not a problem as explained in detail\footnote{See also~\cite[Prop.~6.2]{BDK12} where the generation of non-twisted heterodimensional cycles is proved from the bi-accumulated property.} in~\cite[Sec.~7.1]{BDK12}. Namely, it is proved that since $f$ has a heterodimensional cycle associated with non-real complex leading multipliers, one can always $C^1$-approximate $f$ by  a $C^r$-diffeomorphism $\varphi$  having  a non-twisted heterodimensional cycle associated with a pair of periodic points $P'_\varphi$ and $Q'_\varphi$ with real leading multipliers  which are homoclinically related to the continuations $P_\varphi$ and $Q_\varphi$ of $P$ and $Q$ respectively. Then, applying~\cite[Prop.~6.1]{BDK12} to $\varphi$ we get a diffeomorphism $g$ arbitrarily $C^r$-close to $\varphi$ having the dynamical configuration mentioned above. In particular, note that $g$ has a heterodimensional cycle associated with the saddle $P'_g$ and $Q'_g$ homoclinically related to $P_g$ and $Q_g$ respectively. }

{We explain now following~\cite[Sec.~3.2 and Sec.~6.1.1]{BDK12} how to obtain a $cu$-blender-horseshoe and a heterodimensional cycle as desired after an arbitrarily small $C^r$-perturbation of $g$. First, observe that the dynamical configuration of $g$ is in the assumptions of~\cite[Thm.~3.5]{BDK12}. Hence, from this result, there is a $C^r$-diffeomorphism $h$ arbitrarily $C^r$-close to $g$ with a heterodimensional
cycle associated with the continuation $P'_h$ of $P'_g$ and a $cu$-blender-horseshoe $\Gamma^{cu}_h$ containing a
hyperbolic continuation $S^+_h$ of $S_g$. Moreover, as observed in~\cite[pg.~955]{BDK12}, the dynamical configuration of $g$ implies that the saddle $S^+_h$ can be chosen homoclinically related to $Q'_h$. Thus, by an arbitrarily small $C^r$-pertubation if necessary that we still denote by $h$, one can obtain a heterodimensional cycle associated with $P'_h$ and~$Q'_h$.

Similarly, $h$ also has a $cs$-blender-horseshoe $\Gamma^{cs}_h$  containing a
hyperbolic continuation $S^-_h$ of $S_g$ which can be chosen homoclinically related to $P'_h$.  Unfortunately, this is not notified in~\cite{BDK12} although it is rather folkloric. To justify this, observe that~\cite[Thm.~3.5]{BDK12} is based in~\cite[Thm.~2.4]{BD08}. The main perturbation to prove this result is depicted in~\cite[Fig.~9~and~10]{BD12}. Compare such figures with~\cite[Fig.~2~and~3]{BIR16}. In~\cite[Thm.~2.3]{BIR16}, the authors proved a similar result on the generation of robust cycles from a strong homoclinic intersection following~\cite[Sec.~4]{BD12}. The generation of both $cs$-blender-horseshoe and $cu$-blender-horseshoe is proved in~\cite[Claim~2.2 and~paragraph before Claim~2.3]{BIR16}. %independently of the extra assumption of the existence of a tangencial strong homoclinic intersection.
These blender-horseshoes contain the hyperbolic periodic points $S^-_h$ and $S^+_h$  which are homoclinically related to $P'_h$ and $P'_h$ respectively as mentioned before from the dynamical configuration of $g$. This completes all properties that we need and concludes the proof.}
\end{proof}

%\vspace{0.25cm}
Now, we will conclude the proof of Theorem~\ref{thm1} by proving Proposition~\ref{prop-H123}.

\begin{proof}[Proof of Propostion~\ref{prop-H123}] Let $f$ be a three-dimensional diffeomorphism under the assumpition of Proposition~\ref{prop-H123}. Fix a $C^r$-diffeomorphism $h$ arbitrarily $C^1$-close to $f$ provided by Lemma~\ref{lem-medio}.
Kiriki, Soma and Nakano proved in \cite[Prop.~2.1]{KNS17} that there is a diffeomorphism $g$ arbitrarily $C^r$-close to the diffeomorphism $h$  having a non-transverse equidimensional cycle associated with the continuations $Q'_g$ and $Q_g$ of $Q'_h$ and $Q_h$ respectively\footnote{In order to fit our assumptions to the notation in the referred literature we need to consider $h^{-1}$.}. Thus, $h$ satisfies~\ref{H1}.  But, since $P'_g$ is homoclinically related to the continuation $\Gamma^{cs}_g$ of $\Gamma^{cs}_h$ (which is also a $cs$-blender-horseshoe), we also have~\ref{H2}. To conclude the proof note that $Q'_g$ satisfies~\ref{H3} since $Q$ verifies~\eqref{eq1}.
\end{proof}

To conclude this section we want to remark that the condition~\eqref{eq1} on the multipliers of $Q$ is just imposed to get~\ref{H3}. If this assumption is avoided, we have the following result:

\begin{thm}
Let $f$ be a three-dimensional $C^r$-diffeomorphism exhibiting a
heterodimensional cycle associated  with a pair of hyperbolic
periodic points with non-real complex multipliers. Then $f$ can be $C^1$-approximated by a $C^r$-diffeomorphism $g$  having simultaneously,
\begin{enumerate}
\item \label{tang} a $C^1$-robust homoclinic tangency,
\item \label{het} a $C^1$-robust heterodimensional cycle, and
\item\label{equi} a $C^1$-robust non-transverse equidimensional cycle.
\end{enumerate}
\end{thm}

\begin{proof}
 Under the assumption of this theorem, Proposition~\ref{prop-H123} implies that $f$ can be $C^1$-approximated by  $C^r$-diffeomorphisms satisfying~\ref{H1} and~\ref{H2}. Then, by Proposition~\ref{lem-new}, we can continue to $C^1$-approximate $f$ by a $C^r$-diffeomorphism $h$ with a $C^1$-robust non-transverse equidimensional cycle between a saddle with non-real complex multipliers $Q_h$ and a $cs$-blender-horseshoe $\Gamma_h$. As it is well-known, by using the Inclination lemma, we can $C^r$-approximate $h$ by a diffeomorphism $g$ with a homoclinic tangency associated with a periodic point in $\Gamma_g$. Thus, according to~\cite[Thm.~4.9]{BD12}, there is a diffeomorphism arbitrarily $C^r$-close to $g$ having a $C^1$-robust homoclinic tangency. In particular, this diffeomorphism satisfies simultaneously~\ref{tang} and~\ref{equi}. We will explain that also $g$ satisfies~\ref{het}. To see this, observe that the first step to prove Proposition~\ref{prop-H123} is Lemma~\ref{lem-medio}. This lemma provides a heterodimensional cycle which is $C^1$-robust by the presence of a blender-horseshoe. Thus $h$ and $g$ above satisfy~\ref{het}. This concludes the proof.
\end{proof}
}
%\vspace{-0.1cm}
\subsection{Proof of Theorem~\ref{thm2}} Let $f$ be a
$C^r$-diffeomorphism for $r\geq 2$ with a Tatjer tangency
associated with a dissipative but non-sectional dissipative
periodic point. By a small $C^r$-perturbation, we can assume that,
$f$ is actually  $C^{r+1}$. Now, let us explain in more detail
the results obtained by Tatjer in~\cite{Ta01,GGT07}.

First of all, Tatjer localizes in~\cite[Prop.~{3.5} and 3.7]{Ta01} a sequence
$g_n$ of perturbations of $f$ with a $n$-periodic point $p_n$
having a Bogdanov-Taken bifurcation converging to $f$ as $n$ goes
to infinity. After that, in the proof of~\cite[Prop.~4.1 and
4.5]{Ta01}, Tatjer performs a change of variables around the
$n$-periodic point $p_n$ of $g_n$ in order to calculate a
manageable expression of the return map $h_n=g_n^n$. Since $p_n$
is a {Bogdanov-Takens} bifurcation, $h_n$ has a two-dimensional
invariant center $C^{r}$-manifold $\mathcal{M}_n$ which is
attracting: See~\cite[pg.~293, line 5-7]{Ta01} and \cite[Comments
after Thm.~2]{GGT07}. Actually, in case A, the restriction of
$h_n$ to this manifold is again well approached by the Generalized
H\'enon map (c.f~\cite[Thm.~3]{GST08}). In case B, the
attracting character of the center manifold $\mathcal{M}_n$
follows from the limit return~\cite[Pg.~299, proof of Thm~1,
item~1 and~3]{Ta01} since the surface $y=a+bz+x^2$ is invariant
and every point in $\mathbb{R}^3$ falls by one iteration of limit
map into this surface. The three-dimensional limit return map has
a zero eigenvalue and the study of this family can be reduced to a
family of two-dimensional
endomorphisms~(see~\cite{pumarino2006dynamics} reference therein).
The limit map is obtained by truncating a Taylor expansion in the
manageable expression of the periodic return diffeomorphism. Then,
going from the limit map to the family of return maps, the zero
eigenvalue becomes a real one with a small modulus. Thus, the
invariant center manifold $\mathcal{M}_n$ coming from the
Bogdanov-Takens bifurcation becomes an attracting manifold
providing attracting invariant smooth circle among others after
bifurcation as mentioned in~\cite[pg.~299]{Ta01}. Moreover, by the
results of Broer et al~\cite{broer1996invariant}
(cf.~\cite{Ta01,GGT07}), near a Bogdanov-Takens bifurcation there
{exist} homoclinic tangencies associated with a dissipative saddle
periodic point into the two-dimensional invariant center manifold.
Since the other direction is strong contracting, this periodic
point view in three dimensions is sectional dissipative.
% which allows Tatjer~\cite[pg.~299]{Ta01} to apply the result on the existence
%of strange attractors from~\cite{viana1993strange}.

Similar as it did in \S\ref{sec:proofA}, a three-dimensional historic wandering  domain could be obtained by applying
now~\cite{KS17}. Moreover, notice that, in this case, we
obtain the persistence of homoclinic tangencies directly from the
result of Newhouse~\cite{New79} in dimension 2. This provides a
$C^{r}$-Newhouse domain $\mathcal{N}$ associated with a sectional
dissipative periodic points which are arbitrary $C^r$-close to
$f$ and where maps with non-trivial historic contractive wandering
domains are $C^r$-dense in $\mathcal{N}$. This concludes the proof
of Theorem~\ref{thm2}.

\subsection{Proof of Theorem~\ref{thm3}} \label{sec:proofC}
Because of Theorem~\ref{thm2} we only need to prove the existence
of $C^r$-Newhouse domains of Tatjer homoclinic tangencies
associated with dissipative but non-sectional dissipative periodic
points for $r\geq 1$.

Let us consider a three-dimensional $C^r$-diffeomorphism $f$ having a
non-transverse equidimensional cycle associated with periodic
points $P$ and $Q$ (assumption~\ref{H1}) for $r\geq 2$. Assume
that $Q$ has complex multipliers and all the multipliers of $P$
are real. According to~\cite[Prop.~3.1]{KNS17} such {diffeomorphism}
can be $C^r$-approximated by Tatjer homoclinic tangency
associated with the continuation of the periodic point $P$.
Although~\cite{KNS17} deals with the case that $P$ has unstable
index $2$, the case of unstable index {1} also follows by simply considering $f^{-1}$. However, we cannot conclude, a priori, from
this that  $P$ is a dissipative but non-sectional dissipative
periodic point\footnote{See that this necessary assumption to
apply the results from~\cite{Ta01,GGT07} is missing
in~\cite{KNS17}.}. To do this, we need to impose an extra condition
on the multipliers of $Q$. Namely, we will assume that
\begin{enumerate}[label=(H\arabic*'), start=2]
\item \label{H3'} $Q$ has non-real multipliers and it is dissipative but non-sectional
dissipative.
\end{enumerate}
By Inclination Lemma and a $C^r$-perturbation if necessary, we can
assume that $Q$ has a homoclinic tangency. Hence, according
to~\cite[Sec.~5]{PV94}, one can $C^r$-approximate $f$ by
diffeomorphisms having a homoclinic tangency associated with a
periodic point $P'$ with real multipliers homoclinically related
to {the continuation of} $Q$ and with the same local character, i.e., $P'$ is still dissipative but non-sectional dissipative. Again, by a
$C^r$-perturbation, we can obtain a non-transverse equidimensional
cycle associated with $Q$ and $P'$. Thus, from~\cite{KNS17} we get
now a Tatjer homoclinic tangency associated with the
continuation of $P'$. Summarizing,

\begin{lem}  \label{lem-final}
Let $f$ be  a three-dimensional $C^r$-diffeomorphism ($r\geq 2$)
under the assumptions~\ref{H1} and~\ref{H3'}. Then $f$ can be
$C^r$-approximated by diffeomorphisms displaying a Tatjer
homoclinic tangencies associated with a dissipative but
non-sectional dissipative periodic point.
\end{lem}

To construct a $C^r$-Newhouse domain ($r\geq 1$) we need more
extra assumptions:
\begin{enumerate}[label=(H\arabic*'), start=3]
\item \label{H2'} $P$ is homoclinically related to a blender-horseshoe $\Gamma$.
\end{enumerate}
The above assumptions~\ref{H1}-\ref{H3'}-\ref{H2'} can be
rewritten as follows: either,
\begin{itemize}
\item $f$ satisfies~\ref{H1}-\ref{H3}-\ref{H2} or
\item $f$ satisfies~\ref{H1}, the multipliers of $Q$ are $\lambda$, $\gamma_1$,
$\gamma_2$ with
$$
   |\lambda|<1 < |\gamma| \quad \text{and} \quad   |\lambda \gamma^2| <1
   \quad \text{where} \ \ \lambda\in \mathbb{R}, \ \gamma_{1,2}=\gamma e^{\varphi i} \
    \ \ \varphi\not = 0,\pi
$$
and $P$ is homoclinically related to a $cu$-blender-horseshoe
$\Gamma$.
\end{itemize}
Similar as indicated in the introduction,
assumptions~\ref{H1}-\ref{H3'}-\ref{H2'} implies that $f$ is
$C^r$-approximated by a $C^1$-robust equidimensional tangency
associated with the {continuation} of $Q$ and $\Gamma$
(cf.~\cite[Sec.~4.3]{BD12}). In other words, arbitrarily
$C^r$-close to $f$, there exists a $C^1$-open set  $\mathcal{N}$
of diffeomorphisms such that any $g\in \mathcal{N}$ has a tangency
 between some of the invariant manifolds of  $Q$ and $\Gamma$.
By a $C^r$-perturbation, we can get a smooth diffeomorphism $h$
arbitrarily close to $g$ such that $h$ has a non-transverse
equidimensional cycle associated with $Q$ and some periodic point
$P\in \Gamma$. Thus, we obtain a dense set $\mathcal{D}$ in
$\mathcal{N}$ where any $h\in\mathcal{D}$ is smooth and
satisfies~\ref{H1} and \ref{H3'}. Then, from Lemma~\ref{lem-final}
we can approximate $h$ by Tatjer homoclinic tangencies associated
with dissipative but non-sectional dissipative periodic points.
This proves  the following:

\begin{prop} \label{prop-final}
Arbitrarily $C^r$-close ($r\geq 1$) to a $C^r$-diffeomorphism $f$
satisfying \ref{H1}-\ref{H3'}-\ref{H2'}, there exists a
$C^r$-Newhouse domain $\mathcal{N}$  of  Tatjer homoclinic
tangencies associated with dissipative but non-sectional
dissipative periodic points.
\end{prop}

To complete the proof of Theorem~\ref{thm3} we need to show the
following. Arbitrarily $C^1$-close to a heterodimensional cycle
associated with saddles $P$ and $Q$ with complex multipliers where
$Q$ satisfies~\ref{H3'}, one can find a $C^r$-Newhouse domain of
Tatjer homoclinic tangencies associated with dissipative but
non-sectional dissipative periodic points. As arguing
in~\S\ref{sec:proofA}, heterodimensional cycles in the above
assumptions can be $C^1$-approximated by non-transverse
equidimensional cycles satisfying
assumptions~\ref{H1}-\ref{H3'}-\ref{H2'}. Hence, from
Proposition~\ref{prop-final} we get that $f$ can be
$C^1$-approximated by $C^r$-Newhouse domains of Tatjer homoclinic
tangencies associated with dissipative but non-sectional
dissipative periodic points.

\subsection{Proof of Theorem~\ref{thm4}} \label{sec:proofD}
Let us first recall the result of Colli and Leal
in~\cite{colli1998infinitely} and~\cite{leal2008high}.

\begin{thm}[Colli, Leal] \label{thm-final2} Let $f$ be a $C^\infty$-diffeomorphisms
having a homoclinic tangency associated with a sectional
dissipative periodic point. Then, there exists a $C^\infty$-open
set $\mathcal{U}$ containing $f$ in its closure such that every
$g\in \mathcal{U}$ can be $C^\infty$-approximated by
diffeomorphisms exhibit infinitely many H\'enon-like strange
attractors.
\end{thm}

Now, we will consider a $C^r$-Newhouse domain $\mathcal{N}$
($r\geq 1$) of sectional dissipative periodic points.
%For $r\geq 2$ this set is in the closure of an open set
%of persistence of tangencies (cf.~\cite{PV94,GST93}). For $r=1$, ???
The set  $\mathcal{N}^\infty=\mathcal{N}\cap
\mathrm{Diff}^\infty(M)$ is $C^r$-dense in $\mathcal{N}$ and
$C^\infty$-open in the set of $C^\infty$-diffeomorphisms
$\mathrm{Diff}^\infty(M)$ of the manifold $M$. Moreover, any
diffeomorphism $f\in \mathcal{N}^\infty$ can be
$C^\infty$-approximated by $C^\infty$-diffeomorphisms $g$ having
homoclinic tangencies associated with  sectional dissipative
periodic points.
% This is possible using the Lambda-lemma when we have persistence of tangencies.
Thus, from Theorem~\ref{thm-final2}, $g$ can also be
$C^\infty$-approximated by diffeomorphisms $h$ exhibiting
infinitely many H\'enon-like strange attractors. Consequently,
\emph{there exists a $C^r$-dense set $\mathcal{D}$ in
$\mathcal{N}$ where any $h\in \mathcal{D}$ exhibits infinitely
many H\'enon-like strange attractors}. In what follows, we will
explain how to use this result to conclude Theorem~\ref{thm4}.

First of all, observe that as an immediate consequence of
Theorem~\ref{thm2} and the proof of Theorem~\ref{thm1}, we have
the following:

\begin{prop}
Any $C^r$-Newhouse domain ($r\geq 1$) of type \ref{2} or \ref{3}
in Theorem~\ref{thm4} is also a
%\begin{itemize}
%\item of homoclinc tangencies associated with periodic points
%satisfying~\eqref{eq1},
%\item of Tatjer
%homoclinic tangencies associated with dissipative but
%non-sectional dissipative periodic points
%\end{itemze}
 $C^r$-Newhouse domain of type~\ref{1}.
 \end{prop}

Thus, we can restrict our attention to the case where
$\mathcal{N}$ is a $C^r$-Newhouse domain of type~\ref{1} in
Theorem~\ref{thm4}, i.e., of homoclinic tangencies associated with
sectional dissipative homoclinic periodic points. Recall that, as
notified in the introduction, H\'enon-like strange attractors are
non-trivial attracting homoclinic classes. Since non-trivial
attracting homoclinic classes persist under $C^1$-perturbations,
fixed a finite number $n\in \mathbb{N}$, the above result implies the
existence of an open and dense set $\mathcal{O}_n$ of $\mathcal{N}$
where $n$ different (independent) homoclinic classes of this type
coexists. Taking $\mathcal{R}=\cap \mathcal{O}_n$ we get the
desired residual set and complete the proof.

%\begin{remark}
  {Note that in Theorem~\ref{thm4} we cannot replace non-trivial attracting homoclinic classes by H\'enon like strange attractors. A result asserting the coexistence of infinitely many strange attractors in a residual set of a Newhouse domain is \emph{not expected}. However, the author of this paper proved recently in~\cite{barrientos2021typical} that such result holds in \emph{Berger domains} of persistent homoclinic tangencies to sectional dissipative periodic points. That is, in an open set of parameter families that plays the same role as the notion of Newhouse domain in the free parameter case.}
%\end{}

\subsection*{Acknowledgements}
I thank A.~Raibekas for his unconditional friendship and to
provide the initial idea of this paper. Basically, he wrote the
proof of Theorem~\ref{thm2}. Also, I thank S.~Kiriki, T.~Soma and
Y.~Nakano with whom I started this project about the higher version
of~\cite{KS17} collaborating in the initial version of this paper.
Finally, I thank L.~J.~D\'iaz to introduce the problem on the
coexistence of infinitely many disjoints non-trivial homoclinic
classes and useful conversations to get
Theorem~\ref{thm4}.

\bibliography{wandering-cycles3.bbl}
%\bibliography{bib}

\end{document}